\documentclass[leqno,12pt]{amsart}
\usepackage{amsfonts}
\usepackage{amsmath,amssymb,amsthm}

\newcommand{\R}{\mathbb R}

\renewcommand{\span}{\mathrm{span}}
\newcommand{\tr}{\mathrm{tr}}

\newtheorem{thm}{Theorem}[section]

\newtheorem{prop}[thm]{Proposition}
\theoremstyle{definition}

\theoremstyle{remark}

\newcommand{\ds}{\displaystyle}

\begin{document}

\title[Meridian Surfaces of elliptic or hyperbolic type in the
Minkowski 4-space] {Meridian Surfaces of elliptic or hyperbolic type in the
Four-dimensional Minkowski Space}

\author{Georgi Ganchev and Velichka Milousheva}
\address{Institute of Mathematics and Informatics, Bulgarian Academy of Sciences,
Acad. G. Bonchev Str. bl. 8, 1113 Sofia, Bulgaria}
\email{ganchev@math.bas.bg}
\address{Institute of Mathematics and Informatics, Bulgarian Academy of Sciences,
Acad. G. Bonchev Str. bl. 8, 1113, Sofia, Bulgaria;   "L. Karavelov"
Civil Engineering Higher School, 175 Suhodolska Str., 1373 Sofia,
Bulgaria} \email{vmil@math.bas.bg}


\begin{abstract}
We consider a special class of spacelike surfaces in the Minkowski
4-space which are one-parameter systems of meridians of the rotational
hypersurface with timelike or spacelike axis. We call these  surfaces
meridian surfaces of elliptic or hyperbolic type, respectively.
 On the base of our invariant theory of surfaces we study meridian surfaces with special
invariants and give the complete
classification of  the meridian surfaces with
constant Gauss curvature or constant mean curvature. We also classify the Chen meridian surfaces and the meridian surfaces
with parallel normal bundle.
\end{abstract}

\maketitle

\section{Introduction}

One of the fundamental problems of the contemporary
differential geometry of surfaces and hypersurfaces in the
standard model spaces such as the Euclidean space $\R^n$ and the
pseudo-Euclidean space $\R^n_k$ is the investigation of the basic
invariants characterizing the surfaces. Our aim is to study and
classify various important classes of surfaces in the
four-dimensional Minkowski space $\R^4_1$ characterized by
conditions on their invariants.

An invariant theory of spacelike surfaces in  $\R^4_1$ was
developed by the present authors in \cite {GM5}.
We introduced an invariant linear map $\gamma$ of Weingarten-type  in the tangent plane at
any point of the surface, which generates two invariant functions $k = \det \gamma$ and $\varkappa= -\ds{ \frac{1}{2}}\, \tr \gamma$.
On the base of the map $\gamma$ we  introduced principal lines and a
geometrically determined moving frame field at each point of the
surface. Writing derivative formulas of Frenet-type for this frame
field, we obtained eight invariant functions  $\gamma_1, \,
\gamma_2, \, \nu_1,\, \nu_2, \, \lambda, \, \mu, \, \beta_1,
\beta_2$ and proved a fundamental theorem of Bonnet-type, stating
that these eight invariants under some natural conditions
determine the surface up to a rigid motion in $\R^4_1$.

The basic geometric classes of surfaces in $\R^4_1$   are
characterized by conditions on these invariant functions. For
example, surfaces with flat normal connection are characterized by
the condition $\nu_1 = \nu_2$, minimal surfaces are  described by
$\nu_1 + \nu_2 = 0$, Chen surfaces are characterized by  $\lambda
= 0$, and surfaces with parallel normal bundle are characterized by the condition $\beta_1 = \beta_2 = 0$.

In \cite{GM2} we constructed special two-dimensional  surfaces in the Euclidean 4-space $\R^4$
which are one-parameter systems of meridians of the rotational
hypersurface and called  these surfaces \emph{meridian
surfaces}. We classified the meridian surfaces with
constant Gauss curvature, constant mean curvature, and constant
invariant $k$ \cite{GM2}. In \cite{GM-new} we gave the invariants $\gamma_1, \,
\gamma_2, \, \nu_1,\, \nu_2, \, \lambda, \, \mu, \, \beta_1,
\beta_2$ of the meridian surfaces and on the base of these
invariants we classified  completely the Chen meridian surfaces
and the meridian surfaces with parallel normal bundle.

Similarly to the Euclidean case, in \cite{GM6} we constructed
 two-dimensional  spacelike surfaces in the Minkowski 4-space $\R^4_1$   which are
one-parameter systems of meridians of the rotational hypersurface with
timelike or spacelike axis. We called these surfaces \emph{meridian surfaces of elliptic type} and \emph{meridian surfaces of hyperbolic type},
respectively.
The geometric construction of the meridian surfaces is different from the construction of the standard rotational
surfaces with two-dimensional axis. Hence, the class of meridian surfaces is a new source of examples of two-dimensional surfaces in
$\R^4_1$. In \cite{GM6} we found all marginally trapped  meridian surfaces of elliptic or hyperbolic type.

In \cite{GM7} we continued the study of meridian surfaces in $\R^4_1$ considering a rotational hypersurface with
lightlike  axis and constructed two-dimensional surfaces
which are one-parameter systems of meridians of the rotational
hypersurface. We called these surfaces \emph{meridian surfaces of
parabolic type}.  We calculated their basic  invariants and found all
marginally trapped meridian surfaces of parabolic type.

In the present paper we consider meridian surfaces of elliptic or hyperbolic type in $\R^4_1$  and  calculate the invariants $\gamma_1, \,
\gamma_2, \, \nu_1,\, \nu_2, \, \lambda, \, \mu, \, \beta_1,
\beta_2$ of these surfaces. Using the invariants
we describe and  classify completely  the meridian surfaces  of elliptic or hyperbolic type with
constant Gauss curvature (Theorem \ref{T:Gauss curvature}), with constant mean curvature (Theorem \ref{T:mean curvature}), and  with constant
invariant $k$ (Theorem \ref{T:constant k}). In Theorem \ref{T:Chen}  we classify  the Chen meridian surfaces
and in Theorem \ref{T:parallel}  we give the classification of the meridian surfaces
with parallel normal bundle.

\section{Preliminaries} \label{S:Pre}

Let $\R^4_1$  be the four-dimensional Minkowski space  endowed with the metric
$\langle , \rangle$ of signature $(3,1)$ and  $Oe_1e_2e_3e_4$ be a
fixed orthonormal coordinate system, i.e.  $\langle e_1, e_1 \rangle  =
\langle e_2, e_2 \rangle  = \langle e_3, e_3 \rangle  = 1, \, \langle e_4, e_4 \rangle  = -1$.
A surface $M^2: z = z(u,v), \, \, (u,v) \in {\mathcal D}$
(${\mathcal D} \subset \R^2$) in $\R^4_1$ is said to be
\emph{spacelike} if $\langle , \rangle$ induces  a Riemannian
metric $g$ on $M^2$. Thus at each point $p$ of a spacelike surface
$M^2$ we have the following decomposition:
$$\R^4_1 = T_pM^2 \oplus N_pM^2$$
with the property that the restriction of the metric $\langle ,
\rangle$ onto the tangent space $T_pM^2$ is of signature $(2,0)$,
and the restriction of the metric $\langle , \rangle$ onto the
normal space $N_pM^2$ is of signature $(1,1)$.

Denote by $\nabla'$ and $\nabla$ the Levi Civita connections on $\R^4_1$ and $M^2$, respectively.
Let $x$ and $y$ be vector fields tangent to $M^2$ and $\xi$ be a normal vector field.
The formulas of Gauss and Weingarten give the decompositions of the vector fields $\nabla'_xy$ and
$\nabla'_x \xi$ into tangent and normal components:
$$\begin{array}{l}
\vspace{2mm}
\nabla'_xy = \nabla_xy + \sigma(x,y);\\
\vspace{2mm}
\nabla'_x \xi = - A_{\xi} x + D_x \xi,
\end{array}$$
which define the second fundamental tensor $\sigma$, the normal connection $D$
and the shape operator $A_{\xi}$ with respect to $\xi$.
The mean curvature vector  field $H$ of $M^2$ is defined as
$H = \ds{\frac{1}{2}\,  \tr\, \sigma}$.

Let
$M^2: z=z(u,v), \,\, (u,v) \in \mathcal{D}$ $(\mathcal{D} \subset \R^2)$
be a local parametrization on a
spacelike surface in $\R^4_1$.
The tangent space at an arbitrary point $p=z(u,v)$ of $M^2$ is
$T_pM^2 = \span \{z_u,z_v\}$, where $\langle z_u,z_u \rangle > 0$, $\langle z_v,z_v \rangle > 0$.
We use the standard denotations
$E(u,v)=\langle z_u,z_u \rangle, \; F(u,v)=\langle z_u,z_v
\rangle, \; G(u,v)=\langle z_v,z_v \rangle$ for the coefficients
of the first fundamental form.
Let $\{n_1, n_2\}$  be a normal frame field of $M^2$ such that $\langle
n_1, n_1 \rangle =1$, $\langle n_2, n_2 \rangle = -1$, and the
quadruple $\{z_u,z_v, n_1, n_2\}$ is positively oriented in
$\R^4_1$.
The coefficients of the second fundamental form  $II$ of the surface $M^2$ are given
by the following functions
\begin{equation} \notag
L = \ds{\frac{2}{W}} \left|%
\begin{array}{cc}
\vspace{2mm}
  c_{11}^1 & c_{12}^1 \\
  c_{11}^2 & c_{12}^2 \\
\end{array}%
\right|; \quad
M = \ds{\frac{1}{W}} \left|%
\begin{array}{cc}
\vspace{2mm}
  c_{11}^1 & c_{22}^1 \\
  c_{11}^2 & c_{22}^2 \\
\end{array}%
\right|; \quad
N = \ds{\frac{2}{W}} \left|%
\begin{array}{cc}
\vspace{2mm}
  c_{12}^1 & c_{22}^1 \\
  c_{12}^2 & c_{22}^2 \\
\end{array}%
\right|,
\end{equation}
where
$$\begin{array}{lll}
\vspace{2mm}
c_{11}^1 = \langle z_{uu}, n_1 \rangle; & \qquad   c_{12}^1 = \langle z_{uv},
n_1 \rangle; & \qquad  c_{22}^1 = \langle z_{vv}, n_1 \rangle;\\
\vspace{2mm}
c_{11}^2 = \langle z_{uu}, n_2 \rangle; & \qquad  c_{12}^2 = \langle z_{uv},
n_2 \rangle; & \qquad c_{22}^2 = \langle z_{vv}, n_2 \rangle.
\end{array} $$
The second fundamental
form $II$ is invariant up to the orientation of the tangent space
or the normal space of the surface.

The condition $L = M = N = 0$  characterizes points at which
the space  $\{\sigma(x,y):  x, y \in T_pM^2\}$ is one-dimensional.
We call such points  \emph{flat points} of the surface.
These points are analogous to flat points in the theory of surfaces in $\R^3$ and $\R^4$ \cite{GM2}.
In \cite{GM5} we gave a local geometric description of spacelike surfaces consisting of flat points
 proving that
any spacelike surface consisting of flat points whose mean
curvature vector at any point is a non-zero spacelike vector or
timelike vector  either lies in a hyperplane of $\R^4_1$ or is
part of a developable ruled surface in $\R^4_1$.
Further we consider surfaces free of flat points, i.e.  $(L, M, N) \neq (0,0,0)$.

Using the functions $L$, $M$, $N$ and $E$, $F$, $G$  in \cite{GM5}
we introduced a linear map $\gamma$ of Weingarten type in the
tangent space at any point of $M^2$. The map $\gamma$ is invariant with respect to
changes of parameters on $M^2$ as well as to motions in $\R^4_1$. It
generates two invariant functions
$$k =  \frac{LN - M^2}{EG - F^2}, \qquad
\varkappa = \frac{EN+GL-2FM}{2(EG-F^2)}.$$

It turns out that
the invariant $\varkappa$ is the curvature of the normal
connection of the surface
(see \cite{GM5}).
As in the theory of surfaces in $\R^3$ and $\R^4$  the invariant
$k$ divides the points of $M^2$ into the following  types: \emph{elliptic} ($k > 0$),
\emph{parabolic} ($k = 0$), and \emph{hyperbolic} ($k < 0$).

The second fundamental form $II$ determines conjugate, asymptotic, and
principal tangents at a point $p$ of $M^2$ in the standard way.
A line $c: u=u(q), \; v=v(q); \; q\in J \subset \R$ on $M^2$ is
said to be an \emph{asymptotic line}, respectively a
\textit{principal line}, if its tangent at any point is
asymptotic, respectively  principal. The surface $M^2$ is
parameterized by principal lines if and only if $F=0, \,\, M=0.$

Considering spacelike surfaces in $\R^4_1$ whose mean
curvature vector at any point is a non-zero spacelike vector or
timelike vector, on the base of the principal lines we introduced
a geometrically determined orthonormal frame field $\{x,y,b,l\}$ at each point of such a surface  \cite{GM5}.
The tangent vector fields $x$ and $y$ are collinear with  the principal directions,
the normal vector field $b$ is collinear with the mean curvature vector field $H$.
Writing derivative formulas
of Frenet-type for this frame field, we obtained eight invariant
functions $\gamma_1, \, \gamma_2, \, \nu_1,\, \nu_2, \, \lambda, \, \mu,
\, \beta_1, \beta_2$, which determine the surface up to a rigid motion in $\R^4_1$.

The invariants $\gamma_1, \, \gamma_2, \, \nu_1,\, \nu_2, \, \lambda, \, \mu,
\, \beta_1$, and $\beta_2$ are determined by the geometric frame field $\{x,y,b,l\}$ as follows
\begin{equation}
\begin{array}{l} \label{E:Eq1}
\vspace{2mm}
\nu_1 = \langle \nabla'_xx, b\rangle, \qquad \nu_2 = \langle \nabla'_yy, b\rangle, \qquad \lambda = \langle \nabla'_xy, b\rangle,
\qquad \mu = \langle \nabla'_xy, l\rangle,\\
\vspace{2mm}
\gamma_1 =  \langle \nabla'_xx, y\rangle, \qquad  \gamma_2 = \langle \nabla'_yy, x\rangle, \qquad \beta_1 = \langle \nabla'_xb, l\rangle, \qquad
\beta_2 = \langle \nabla'_yb, l\rangle.
\end{array}
\end{equation}

The invariants $k$, $\varkappa$, and the Gauss curvature $K$ of
$M^2$ are expressed by the functions $\nu_1, \nu_2, \lambda, \mu$
as follows:
\begin{equation} \notag
k = - 4\nu_1\,\nu_2\,\mu^2, \quad \quad \varkappa = (\nu_1-\nu_2)\mu, \quad \quad
K = \varepsilon (\nu_1\,\nu_2- \lambda^2 + \mu^2),
\end{equation}
where $\varepsilon = sign \langle H, H \rangle$.
The  norm $\Vert H \Vert$ of the mean curvature
vector is expressed as
\begin{equation} \notag
\Vert H \Vert = \displaystyle{ \frac{|\nu_1 + \nu_2|}{2} = \frac{\sqrt{\varkappa^2-k}}{2 |\mu |}}.
\end{equation}
If $M^2$ is a spacelike surface whose mean curvature vector at any point is a non-zero spacelike vector or timelike vector, then
$M^2$ is minimal if and only if $\nu_1 + \nu_2 = 0$.

The geometric meaning of the invariant $\lambda$ is connected with the notion of Chen submanifolds.
Let $M$ be an $n$-dimensional submanifold of
$(n+m)$-dimensional Riemannian manifold $\widetilde{M}$ and $\xi$
be a normal vector field of $M$. B.-Y. Chen \cite{Chen1} defined
 the \emph{allied vector field} $a(\xi)$ of $\xi$  by the
formula
$$a(\xi) = \ds{\frac{\|\xi\|}{n} \sum_{k=2}^m \{\tr(A_1 A_k)\}\xi_k},$$
where $\{\xi_1 = \ds{\frac{\xi}{\|\xi\|}},\xi_2, \dots,  \xi_m \}$ is an
orthonormal base of the normal space of $M$, and $A_i = A_{\xi_i},
\,\, i = 1,\dots, m$ is the shape operator with respect to
$\xi_i$. The allied vector field $a(H)$ of the mean
curvature vector field $H$ is called the \emph{allied mean
curvature vector field} of $M$ in $\widetilde{M}$. B.-Y. Chen
defined  the $\mathcal{A}$-submanifolds to be those submanifolds
of $\widetilde{M}$ for which
 $a(H)$ vanishes identically \cite{Chen1}.
In \cite{GVV1}, \cite{GVV2} the $\mathcal{A}$-submanifolds are
called \emph{Chen submanifolds}. It is easy to see that minimal
submanifolds, pseudo-umbilical submanifolds and hypersurfaces are
Chen submanifolds. These Chen submanifolds are said to be trivial
Chen-submanifolds. In \cite{GM5} we showed that if  $M^2$ is a spacelike surface
in $\R^4_1$ with spacelike or timelike mean curvature vector field then the allied mean curvature vector field of $M^2$
is
$$a(H) = \ds{\frac{\sqrt{\varkappa^2-k}}{2} \,\lambda \, l}.$$
Hence, if $M^2$ is free of minimal points, then $a(H) = 0$ if and
only if $\lambda = 0$. This gives the geometric meaning of the
invariant $\lambda$:  $M^2$ is a non-trivial  Chen
surface if and only if the invariant $\lambda$ is zero.

Now we shall discuss the geometric meaning of the invariants $\beta_1$ and $\beta_2$. It follows from \eqref{E:Eq1} that
\begin{equation} \notag
\begin{array}{ll}
\vspace{2mm}
\nabla'_xb = - \nu_1\,x - \lambda\,y - \beta_1\,l;  & \qquad
\nabla'_xl = - \mu\,y - \beta_1\,b;\\
\vspace{2mm}
\nabla'_yb = - \lambda\,x - \nu_2\,y - \beta_2\,l; & \qquad \nabla'_yl = - \mu\,x - \beta_2\,b.
\end{array}
\end{equation}
Hence, $\beta_1 = \beta_2 = 0$ if and only if $D_xb = D_yb = 0$ (or equivalently, $D_xl = D_yl = 0$).

A normal vector field $\xi$ is said to be \emph{parallel in the normal bundle} (or simply \emph{parallel}) \cite{Chen2},
if $D_x\xi = 0$ holds identically for any tangent vector field $x$.
Hence, the invariants  $\beta_1$ and $\beta_2$  are identically  zero if and only if the geometric normal vector fields $b$ and $l$ are parallel in the normal bundle.

Surfaces  admitting a geometric normal frame field $\{b, l\}$ of parallel normal vector fields, we shall call \emph{surfaces with parallel normal bundle}.
They are characterized by the condition $\beta_1 = \beta_2 = 0$.
Note that if $M^2$ is a surface  free of minimal points with parallel mean curvature vector field (i.e. $DH = 0$),
then $M^2$ is a surface with parallel normal bundle, but the converse is not true in general.
It is true only in the case $\Vert H \Vert = const$.

\section{Invariants of meridian surfaces of elliptic or hyperbolic type}

In \cite{GM2} we constructed a family of surfaces lying on a
standard rotational hypersurface in the four-dimensional Euclidean
space $\R^4$. These surfaces are one-parameter systems of
meridians of the rotational hypersurface, that is why we called
them \emph{meridian surfaces}. In \cite{GM6} we used the idea from the Euclidean case to construct
special families of two-dimensional  spacelike surfaces lying
 on rotational hypersurfaces in $\R^4_1$  with
timelike or spacelike axis. The construction was the following.

Let $f = f(u), \,\, g = g(u)$ be smooth functions, defined in an
interval $I \subset \R$, such that $\dot{f}^2(u) - \dot{g}^2(u) > 0,
\,\, u \in I$. We assume that $f(u)>0, \,\, u \in I$. The standard
rotational hypersurface $\mathcal{M}'$ in $\R^4_1$, obtained by
the rotation of the meridian curve $m: u \rightarrow (f(u), g(u))$
about the $Oe_4$-axis,  is parameterized as follows:
$$\mathcal{M}': Z(u,w^1,w^2) = f(u)\, \cos w^1 \cos w^2 \,e_1 +  f(u)\, \cos w^1 \sin w^2 \,e_2 + f(u)\, \sin w^1 \,e_3 + g(u) \,e_4.$$
The rotational hypersurface $\mathcal{M}'$ is a two-parameter system of meridians. Let $w^1 =
w^1(v)$, $w^2=w^2(v), \,\, v \in J, \, J \subset \R$. We consider the two-dimensional
 surface $\mathcal{M}'_m$ lying on $\mathcal{M}'$, constructed  in the following way:
\begin{equation}  \notag
\mathcal{M}'_m: z(u,v) = Z(u,w^1(v),w^2(v)), \quad u \in I, \, v \in J.
\end{equation}
$\mathcal{M}'_m$ is a
one-parameter system of meridians of $\mathcal{M}'$. We call
$\mathcal{M}'_m$ a  \emph{meridian surface of elliptic type}.

If we denote $l(w^1,w^2) = \cos w^1 \cos w^2 \,e_1 + \cos w^1 \sin w^2 \,e_2 + \sin w^1 \,e_3$, then the surface $\mathcal{M}'_m$ is parameterized by
\begin{equation} \label{E:Eq2}
\mathcal{M}'_m: z(u,v) = f(u) \, l(v) + g(u)\, e_4, \quad u \in I, \, v \in J.
\end{equation}

Note that $l(w^1,w^2)$ is  the unit position vector of the 2-dimensional sphere $S^2(1)$ lying in the Euclidean space $\R^3 = \span \{e_1, e_2, e_3\}$
and centered at the origin $O$.

\vskip 2mm
In a similar way we consider meridian surfaces lying on the
 rotational hypersurface in $\R^4_1$ with spacelike axis.
Let $f = f(u), \,\, g = g(u)$ be smooth functions, defined in an
interval $I \subset \R$, such that $\dot{f}^2(u) + \dot{g}^2(u)
>0$, $f(u)>0, \,\, u \in I$.
The rotational hypersurface $\mathcal{M}''$ in $\R^4_1$, obtained
by the rotation of the meridian curve $m: u \rightarrow (f(u),
g(u))$ about the $Oe_1$-axis  is parameterized as follows:
$$\mathcal{M}'': Z(u,w^1,w^2) = g(u) \,e_1 +  f(u)\, \cosh w^1 \cos w^2 \,e_2 +  f(u)\, \cosh w^1 \sin w^2 \,e_3+ f(u)\, \sinh w^1 \,e_4.$$
If $w^1 = w^1(v), \, w^2=w^2(v), \,\, v \in J, \,J \subset
\R$, we construct a surface $\mathcal{M}''_m$ in $\R^4_1$ in
the following way:
\begin{equation}  \notag
\mathcal{M}''_m: z(u,v) = Z(u,w^1(v),w^2(v)),\quad u \in I, \, v \in J.
\end{equation}
 $\mathcal{M}''_m$ is a one-parameter system of
meridians of $\mathcal{M}''$. We call $\mathcal{M}''_m$ a
\emph{meridian surfaces of hyperbolic type}.

If we denote $l(w^1,w^2) = \cosh w^1 \cos w^2 \,e_2 +  \cosh w^1 \sin w^2 \,e_3 + \sinh w^1 \,e_4$, then the surface $\mathcal{M}''_m$
is given by
\begin{equation} \label{E:Eq3}
\mathcal{M}''_m: z(u,v) = f(u) \, l(v) + g(u)\, e_1, \quad u \in I, \, v \in J,
\end{equation}
$l(w^1,w^2)$ being the unit position vector of the  de Sitter space  $S^2_1(1)$  in the Minkowski space $\R^3_1 = \span \{e_2, e_3, e_4\}$,
i.e. $S^2_1(1) = \{ V \in \R^3_1: \langle V, V \rangle = 1\}$.

In \cite{GM6} we found all  marginally trapped meridian surfaces of elliptic or hyperbolic type.
In the present section we shall find the geometric invariant
functions $\gamma_1, \, \gamma_2, \, \nu_1,\, \nu_2, \, \lambda, \, \mu, \, \beta_1, \beta_2$ of the meridian surfaces of elliptic or hyperbolic type.

\vskip 2mm
\emph{Elliptic case:}

First we consider the surface $\mathcal{M}'_m$  parameterized by \eqref{E:Eq2}.
We assume that the smooth curve $c: l = l(v) = l(w^1(v),w^2(v)), \, v \in J$  on
$S^2(1)$ is parameterized by the arc-length, i.e. $\langle l'(v), l'(v) \rangle = 1$. Let
 $t(v) = l'(v)$  be the tangent vector field of $c$. Since $\langle t(v), t(v) \rangle = 1$,
 $\langle l(v), l(v) \rangle = 1$, and $\langle t(v), l(v) \rangle = 0$,
 there exists a unique (up to a sign)
 vector field $n(v)$, such that
$\{ l(v), t(v), n(v)\}$ is an orthonormal frame field in $\R^3$. With respect to this
 frame field we have the following Frenet formulas of $c$ on $S^2(1)$:
\begin{equation} \label{E:Eq4}
\begin{array}{l}
\vspace{2mm}
l' = t;\\
\vspace{2mm}
t' = \kappa \,n - l;\\
\vspace{2mm} n' = - \kappa \,t,
\end{array}
\end{equation}
 where $\kappa (v)= \langle t'(v), n(v) \rangle$ is the spherical curvature of $c$.

Without loss of generality we assume that $\dot{f}^2(u) - \dot{g}^2(u) = 1$.
The tangent space of $\mathcal{M}'_m$ is spanned by the vector fields:
$$z_u = \dot{f} \,l + \dot{g}\,e_4; \qquad  z_v = f\,t,$$
so, the coefficients of the first fundamental form of $\mathcal{M}'_m$  are $E = 1; \, F = 0; \, G = f^2(u) >0$.
Hence, the first fundamental form is positive
definite, i.e. $\mathcal{M}'_m$ is a spacelike surface.

Denote  $X = z_u,\,\, Y = \ds{\frac{z_v}{f} = t}$ and
consider the following orthonormal normal frame field:
$$n_1 = n(v); \qquad n_2 = \dot{g}(u)\,l(v) + \dot{f}(u) \, e_4.$$
Thus we obtain a frame field $\{X,Y, n_1, n_2\}$ of $\mathcal{M}'_m$, such that $\langle n_1, n_1 \rangle =1$,
$\langle n_2, n_2 \rangle =- 1$, $\langle n_1, n_2 \rangle =0$.

Taking into account \eqref{E:Eq4} we get
the following derivative formulas:
\begin{equation} \label{E:Eq5}
\begin{array}{ll}
\vspace{2mm} \nabla'_XX = \qquad \qquad \qquad \qquad
\kappa_m\,n_2; & \qquad
\nabla'_X n_1 = 0;\\
\vspace{2mm} \nabla'_XY = 0;  & \qquad
\nabla'_Y n_1 = \ds{\quad \quad \, - \frac{\kappa}{f}\,Y};\\
\vspace{2mm} \nabla'_YX = \quad\quad
\quad\ds{\frac{\dot{f}}{f}}\,Y;  & \qquad
\nabla'_X n_2 = \kappa_m \,X;\\
\vspace{2mm} \nabla'_YY = \ds{- \frac{\dot{f}}{f}\,X \quad\quad +
\frac{\kappa}{f}\,n_1 + \frac{\dot{g}}{f} \, n_2}; & \qquad
\nabla'_Y n_2 = \ds{ \quad \quad \quad  \frac{\dot{g}}{f}\,Y},
\end{array}
\end{equation}
where $\kappa_m$ denotes the curvature of the
meridian curve $m$, i.e.
$\kappa_m (u)= \dot{f}(u) \ddot{g}(u) - \dot{g}(u) \ddot{f}(u)$.

The invariants $k$, $\varkappa$, and the Gauss curvature $K$ are given by the following formulas \cite{GM6}:
$$k = - \frac{\kappa_m^2(u) \, \kappa^2(v)}{f^2(u)}; \qquad \varkappa = 0;
\qquad K = - \frac{\ddot{f}(u)}{f(u)}.$$

The equality $\varkappa = 0$ implies the following statement.

\begin{prop}
The meridian surface of elliptic type $\mathcal{M}'_m$, defined by \eqref{E:Eq2}, is a surface with
flat normal connection.
\end{prop}

We distinguish the following  three cases:
\vskip 2mm
I. $\kappa(v) = 0$, i.e. the curve $c$ is a great circle on $S^2(1)$. In
this case $n_1 = const$, and $\mathcal{M}'_m$ is a planar surface lying in
the constant 3-dimensional space spanned by $\{X, Y, n_2\}$.

\vskip 2mm
II. $\kappa_m(u) = 0$, i.e. the meridian curve $m$ is part of
a straight line. In such case $k = \varkappa = K = 0$, and $\mathcal{M}'_m$
is a developable ruled surface.

\vskip 2mm
III. $\kappa_m(u) \, \kappa(v) \neq 0$, i.e. $c$ is not a
great circle on $S^2(1)$ and $m$ is not a straight line.

\vskip 2mm
In the first two cases the surface $\mathcal{M}'_m$ consists of flat points. So, we consider the  third (general) case,
 i.e. we assume that $\kappa_m \neq 0$ and $\kappa \neq 0$.

It follows from \eqref{E:Eq5} that the mean curvature vector field $H$ of $\mathcal{M}'_m$ is expressed as
$$H = \frac{\kappa}{2f}\, n_1 + \frac{\dot{g} + f \kappa_m}{2f} \, n_2.$$
Using that  $\dot{g}^2(u) = \dot{f}^2(u) - 1$  and $\kappa_m (u) = \ds{\frac{\ddot{f}(u)}{\sqrt{\dot{f}^2(u) - 1}}}$, we get
\begin{equation} \label{E:Eq6}
H = \frac{\kappa}{2f}\, n_1 + \frac{ f \ddot{f} + \dot{f}^2 - 1}{2f \sqrt{\dot{f}^2 - 1}} \, n_2.
\end{equation}

Since $\kappa \neq 0$ the surface $\mathcal{M}'_m$ is non-minimal, i.e. $H \neq 0$.
The  case $\mathcal{M}'_m$ is a marginally trapped  surface, i.e. $H \neq 0$ and $\langle H, H \rangle = 0$ is described in \cite{GM6}.
So, here we consider the case  $\langle H, H \rangle \neq 0$.

Note that the orthonormal frame field $\{X,Y, n_1, n_2\}$ of $\mathcal{M}'_m$ is not the  geometric frame field  defined in Section \ref{S:Pre}.
The principal tangents of $\mathcal{M}'_m$ are
\begin{equation}\notag
x = \ds{\frac{X+Y}{\sqrt{2}}}; \qquad   y = \ds{\frac{- X + Y}{\sqrt{2}}}.
\end{equation}
In the case $\langle H, H \rangle > 0$, i.e. $\kappa^2 (\dot{f}^2 - 1) - (f \ddot{f} + \dot{f}^2 - 1)^2 > 0$, the geometric normal frame field  $\{b,l\}$ is given by
\begin{equation}\notag
\begin{array}{l}
\vspace{2mm}
b = \ds{\frac{1}{\sqrt{\kappa^2 (\dot{f}^2 - 1) - (f \ddot{f} + \dot{f}^2 - 1)^2}} \left( \kappa \sqrt{\dot{f}^2 - 1}\, n_1 + (f \ddot{f} + \dot{f}^2 - 1)\,n_2 \right)}; \\
\vspace{2mm}
l = \ds{\frac{1}{\sqrt{\kappa^2 (\dot{f}^2 - 1) - (f \ddot{f} + \dot{f}^2 - 1)^2}} \left( (f \ddot{f} + \dot{f}^2 - 1) \, n_1 + \kappa \sqrt{\dot{f}^2 - 1}\,n_2 \right)}.
\end{array}
\end{equation}
In this case the normal vector fields $b$ and $l$ satisfy $\langle b, b \rangle = 1$, $\langle b, l \rangle = 0$, $\langle l, l \rangle = -1$.

In the case $\langle H, H \rangle < 0$,  i.e. $\kappa^2 (\dot{f}^2 - 1) - (f \ddot{f} + \dot{f}^2 - 1)^2 < 0$,  the geometric normal frame field  $\{b,l\}$ is given by
\begin{equation}\notag
\begin{array}{l}
\vspace{2mm}
b = \ds{- \frac{1}{\sqrt{(f \ddot{f} + \dot{f}^2 - 1)^2 - \kappa^2 (\dot{f}^2 - 1)}} \left( \kappa \sqrt{\dot{f}^2 - 1}\, n_1 + (f \ddot{f} + \dot{f}^2 - 1)\,n_2 \right)}; \\
\vspace{2mm}
l = \ds{\frac{1}{\sqrt{(f \ddot{f} + \dot{f}^2 - 1)^2 - \kappa^2 (\dot{f}^2 - 1)}} \left( (f \ddot{f} + \dot{f}^2 - 1) \, n_1 + \kappa \sqrt{\dot{f}^2 - 1}\,n_2 \right)}.
\end{array}
\end{equation}
In this case we have $\langle b, b \rangle = - 1$, $\langle b, l \rangle = 0$, $\langle l, l \rangle = 1$.

Applying formulas \eqref{E:Eq1} for the geometric frame field $\{x,y, b, l\}$ of $\mathcal{M}'_m$ and derivative formulas \eqref{E:Eq5},
we obtain the following  invariants of  $\mathcal{M}'_m$:
\begin{equation} \label{E:Eq7}
\begin{array}{l}
\vspace{2mm}
\gamma_1 = \gamma_2 = \ds{- \frac{\dot{f}}{\sqrt{2}f}};\\
\vspace{2mm}
\nu_1  = \nu_2 = \ds{ \frac{1}{2f \sqrt{\dot{f}^2 - 1}} \sqrt{ \varepsilon(\kappa^2 (\dot{f}^2 - 1) - (f \ddot{f} + \dot{f}^2 - 1)^2)}};\\
\vspace{2mm}
\lambda = \ds{ \varepsilon \frac{\kappa^2 (\dot{f}^2 - 1) + f^2 \ddot{f}^2 - (\dot{f}^2 - 1)^2}{2f \sqrt{\dot{f}^2 - 1} \sqrt{\varepsilon (\kappa^2 (\dot{f}^2 - 1) - (f \ddot{f} + \dot{f}^2 - 1)^2)}}}; \\
\vspace{2mm}
\mu = \ds{\frac{\kappa \ddot{f}}{\sqrt{\varepsilon (\kappa^2 (\dot{f}^2 - 1) - (f \ddot{f} + \dot{f}^2 - 1)^2)}}}; \\
\vspace{2mm}
\beta_1 = \ds{\frac{ - (\dot{f}^2 - 1)}{\sqrt{2} \varepsilon(\kappa^2 (\dot{f}^2 - 1) - (f \ddot{f} + \dot{f}^2 - 1)^2)}} \left( \kappa\, \frac{d}{du}\left( \frac{f \ddot{f} + \dot{f}^2 - 1}{\sqrt{\dot{f}^2 - 1}} \right) -  \frac{d}{dv}(\kappa) \, \frac{f \ddot{f} + \dot{f}^2 - 1}{f \sqrt{\dot{f}^2 - 1}} \right); \\
\vspace{2mm}
\beta_2 = \ds{\frac{(\dot{f}^2 - 1)}{\sqrt{2} \varepsilon(\kappa^2 (\dot{f}^2 - 1) - (f \ddot{f} + \dot{f}^2 - 1)^2)}} \left( \kappa\, \frac{d}{du}\left( \frac{f \ddot{f} + \dot{f}^2 - 1}{\sqrt{\dot{f}^2 - 1}} \right) + \frac{d}{dv}(\kappa) \, \frac{f \ddot{f} + \dot{f}^2 - 1}{f \sqrt{\dot{f}^2 - 1}} \right),
\end{array}
\end{equation}
where $\varepsilon = sign \langle H, H \rangle$.

\vskip 3mm
\emph{Hyperbolic case:}

Let   $\mathcal{M}''_m$  be the surface parameterized by \eqref{E:Eq3}.
Assume that the curve $c: l = l(v) = l(w^1(v),w^2(v)), \, v \in J$  on
$S^2_1(1)$ is parameterized by the arc-length, i.e. $\langle l'(v), l'(v) \rangle = 1$.
Similarly to the previous case we consider an  orthonormal frame field
$\{ l(v), t(v), n(v)\}$ in $\R^3_1$, such that $t(v) = l'(v)$ and $\langle n(v), n(v) \rangle = -1$.
With respect to this  frame field we have the following decompositions of the vector fields $l'(v)$, $t'(v)$, $n'(v)$:
\begin{equation} \label{E:Eq8}
\begin{array}{l}
\vspace{2mm}
l' = t;\\
\vspace{2mm}
t' = - \kappa \,n - l;\\
\vspace{2mm} n' = - \kappa \,t,
\end{array}
\end{equation}
which can be considered as Frenet formulas of $c$ on $S^2_1(1)$.
The function $\kappa (v)= \langle t'(v), n(v) \rangle$ is the spherical
curvature of $c$ on $S^2_1(1)$.

We assume that $\dot{f}^2(u) + \dot{g}^2(u) = 1$.
Denote  $X = z_u = \dot{f} \,l + \dot{g}\,e_1, \,\, Y = \ds{\frac{z_v}{f} = t}$ and
consider the orthonormal normal frame field defined by:
$$n_1 = \dot{g}(u)\,l(v) - \dot{f}(u) \, e_1; \qquad n_2 =  n(v).$$
Thus we obtain a frame field $\{X,Y, n_1, n_2\}$ of $\mathcal{M}''_m$,
such that $\langle n_1, n_1 \rangle =1$, $\langle n_2, n_2 \rangle =
- 1$, $\langle n_1, n_2 \rangle =0$.

Taking into account \eqref{E:Eq8} we get
the following derivative formulas:
\begin{equation} \label{E:Eq9}
\begin{array}{ll}
\vspace{2mm} \nabla'_XX = \qquad \qquad \quad - \kappa_m\,n_1; & \qquad
\nabla'_X n_1 = \kappa_m \,X;\\
\vspace{2mm} \nabla'_XY = 0;  & \qquad
\nabla'_Y n_1 = \ds{\quad \quad \;\; \frac{\dot{g}}{f}\,Y};\\
\vspace{2mm} \nabla'_YX = \quad\quad
\quad\ds{\frac{\dot{f}}{f}}\,Y;  & \qquad
\nabla'_X n_2 = 0;\\
\vspace{2mm} \nabla'_YY = \ds{- \frac{\dot{f}}{f}\,X \quad\quad -
\frac{\dot{g}}{f}\,n_1 - \frac{\kappa}{f} \, n_2}; & \qquad
\nabla'_Y n_2 = \ds{ \quad \quad \quad  - \frac{\kappa}{f}\,Y},
\end{array}
\end{equation}
where $\kappa_m$ is the curvature of the meridian curve $m$.

The invariants $k$, $\varkappa$, and the Gauss curvature $K$ of the meridian surface  $\mathcal{M}''_m$ are expressed
by the curvatures $\kappa_m(u)$, $\kappa (v)$, and the function $f(u)$ in the same way as the invariants of the meridian surface of elliptic type, i.e.
$$k = - \frac{\kappa_m^2(u) \, \kappa^2(v)}{f^2(u)}; \qquad \varkappa = 0; \qquad K = - \frac{\ddot{f}(u)}{f(u)}.$$

The following statement holds, since $\varkappa = 0$.

\begin{prop}
The meridian surface of hyperbolic type $\mathcal{M}''_m$, defined by \eqref{E:Eq3}, is a surface with
flat normal connection.
\end{prop}

Again we have the following  three cases:
\vskip 2mm
I. $\kappa(v) = 0$. In this case $n_2 = const$, and $\mathcal{M}''_m$ is a planar surface lying in
the constant 3-dimensional space spanned by $\{X, Y, n_1\}$.

\vskip 2mm
II. $\kappa_m(u) = 0$. In such case $k = \varkappa = K = 0$, and $\mathcal{M}''_m$
is a developable ruled surface.

\vskip 2mm
III. $\kappa_m(u) \, \kappa(v) \neq 0$.

\vskip 2mm
In the first two cases  $\mathcal{M}''_m$ is a surface consisting of flat points. So, we consider the  third (general) case,
 i.e. we assume that $\kappa_m \neq 0$ and $\kappa \neq 0$.

Using \eqref{E:Eq9} we get that the mean curvature vector field $H$ of $\mathcal{M}''_m$ is
$$H = - \frac{\dot{g} + f \kappa_m}{2f}\, n_1 -  \frac{\kappa}{2f}\, n_2.$$
Having in mind that  $\dot{g}^2(u) = 1- \dot{f}^2(u)$  and $\kappa_m (u) = \ds{- \frac{\ddot{f}(u)}{\sqrt{1 - \dot{f}^2(u)}}}$, we obtain
\begin{equation} \label{E:Eq10}
H = \frac{ f \ddot{f} + \dot{f}^2 - 1}{2f \sqrt{1 - \dot{f}^2}}\, n_1 -  \frac{\kappa}{2f}\, n_2.
\end{equation}

The surface $\mathcal{M}''_m$ is non-minimal, since $\kappa \neq 0$.
The  case $\mathcal{M}''_m$ is  marginally trapped is described in \cite{GM6}.
So, we consider the case  $\langle H, H \rangle \neq 0$, i.e. $(f \ddot{f} + \dot{f}^2 - 1)^2 - \kappa^2 (1 - \dot{f}^2) \neq 0$.

Similarly to the elliptic case, we find the geometric frame field $\{x,y, b, l\}$ of $\mathcal{M}''_m$.
Applying formulas \eqref{E:Eq1} for this  frame field  and using derivative formulas \eqref{E:Eq9},
we obtain the following  invariants of  $\mathcal{M}''_m$:

\begin{equation} \label{E:Eq11}
\begin{array}{l}
\vspace{2mm}
\gamma_1 = \gamma_2 = \ds{- \frac{\dot{f}}{\sqrt{2}f}};\\
\vspace{2mm}
\nu_1  = \nu_2 = \ds{ \frac{1}{2f \sqrt{1 - \dot{f}^2}} \sqrt{ \varepsilon( (f \ddot{f} + \dot{f}^2 - 1)^2 - \kappa^2 (1 - \dot{f}^2))}};\\
\vspace{2mm}
\lambda = \ds{ \varepsilon \frac{- \kappa^2 (1 - \dot{f}^2) - f^2 \ddot{f}^2 + (1 - \dot{f}^2)^2}{2f \sqrt{1 -\dot{f}^2} \sqrt{\varepsilon ( (f \ddot{f} + \dot{f}^2 - 1)^2 - \kappa^2 (1 - \dot{f}^2))}}}; \\
\vspace{2mm}
\mu = \ds{\frac{\kappa \ddot{f}}{\sqrt{\varepsilon ( (f \ddot{f} + \dot{f}^2 - 1)^2 - \kappa^2 (1 - \dot{f}^2))}}}; \\
\vspace{2mm}
\beta_1 = \ds{\frac{ - (1 - \dot{f}^2)}{\sqrt{2} \varepsilon( (f \ddot{f} + \dot{f}^2 - 1)^2 - \kappa^2 (1 - \dot{f}^2))}} \left( \kappa\, \frac{d}{du}\left( \frac{f \ddot{f} + \dot{f}^2 - 1}{\sqrt{1 - \dot{f}^2}} \right) -  \frac{d}{dv}(\kappa) \, \frac{f \ddot{f} + \dot{f}^2 - 1}{f \sqrt{1 - \dot{f}^2}} \right); \\
\vspace{2mm}
\beta_2 = \ds{\frac{(1 - \dot{f}^2)}{\sqrt{2} \varepsilon((f \ddot{f} + \dot{f}^2 - 1)^2 - \kappa^2 (1 - \dot{f}^2))}} \left( \kappa\, \frac{d}{du}\left( \frac{f \ddot{f} + \dot{f}^2 - 1}{\sqrt{1 - \dot{f}^2}} \right) + \frac{d}{dv}(\kappa) \, \frac{f \ddot{f} + \dot{f}^2 - 1}{f \sqrt{1 - \dot{f}^2}} \right),
\end{array}
\end{equation}
where $\varepsilon = sign \langle H, H \rangle$.

\vskip 2mm
In the following sections, using the invariants of the meridian surfaces $\mathcal{M}'_m$ and $\mathcal{M}''_m$,  we shall describe and classify some special classes of meridian surfaces
of elliptic or hyperbolic type.

\section{Meridian surfaces with constant Gauss curvature}

The study of surfaces with constant Gauss curvature is one of the
main topics in classical differential geometry. Surfaces
with constant Gauss curvature in Minkowski space have drawn the interest of many
geometers, see for example \cite{GalMarMil},  \cite{Lop}, and the references therein.

In the present section we give the classification of the meridian surfaces of elliptic or hyperbolic type in $\R^4_1$ with constant Gauss curvature.

Let $\mathcal{M}'_m$ and  $\mathcal{M}''_m$ be  meridian surfaces of elliptic and hyperbolic type, respectively.
The Gauss curvature in both cases depends only on the meridian curve $m$ and is expressed by the formula
\begin{equation} \label{E:Eq12}
K = - \frac{\ddot{f}(u)}{f(u)}.
\end{equation}

\begin{thm} \label{T:Gauss curvature}
Let $\mathcal{M}'_m$ (resp.  $\mathcal{M}''_m$) be a meridian surface of elliptic (resp. hyperbolic) type from the general class.
Then $\mathcal{M}'_m$ (resp.  $\mathcal{M}''_m$) has constant non-zero Gauss curvature $K$ if and only
if the meridian $m$ is given by
$$\begin{array}{ll}
\vspace{2mm}
f(u) = \alpha \cos \sqrt{K} u + \beta \sin \sqrt{K} u, & \textrm{if} \quad K >0;\\
\vspace{2mm} f(u) = \alpha \cosh \sqrt{-K} u + \beta \sinh
\sqrt{-K} u, & \textrm{if} \quad K <0,
\end{array}$$
where $\alpha$ and $\beta$ are constants, $g(u)$ is defined by $\dot{g}(u) = \sqrt{\dot{f}^2(u) - 1}$ in the elliptic case and
$g(u)$ is defined by $\dot{g}(u) = \sqrt{1 - \dot{f}^2(u)}$ in the hyperbolic case.
\end{thm}

\noindent {\it Proof:} Using \eqref{E:Eq12}  we
obtain that the Gauss curvature $K = const \neq 0$  if and
only if the function $f(u)$ satisfies the following differential
equation
$$\ddot{f}(u) + K f(u) = 0.$$
The general solution of the above equation is given by
$$\begin{array}{ll}
\vspace{2mm}
f(u) = \alpha \cos \sqrt{K} u + \beta \sin \sqrt{K} u, & \quad \textrm{if} \quad K >0;\\
\vspace{2mm} f(u) = \alpha \cosh \sqrt{-K} u + \beta \sinh
\sqrt{-K} u, & \quad \textrm{if}  \quad K <0,
\end{array}$$
where $\alpha$ and $\beta$ are constants.
In the case of meridian surface of elliptic type the function $g(u)$ is
determined by $\dot{g}(u) = \sqrt{\dot{f}^2(u) - 1}$ and in the case of meridian surface of hyperbolic type $\dot{g}(u) = \sqrt{1 - \dot{f}^2(u)}$.

\qed

\section{Meridian surfaces with constant mean curvature}

Constant mean curvature surfaces in arbitrary spacetime
are important objects for their special role in the theory of
general relativity. The study of constant mean curvature surfaces
(CMC surfaces) involves not only geometric methods but also PDE
and complex analysis, that is why the theory of CMC surfaces is of
great interest not only for mathematicians but also for physicists
and engineers. Surfaces with constant mean curvature in
Minkowski space have been studied intensively in  the last years. See for example \cite{Liu-Liu-1}, \cite{Lop-2}, \cite{Sa}, \cite{Chav-Can}, \cite{Bran}.

In this section we classify the  meridian surfaces of elliptic or hyperbolic type with constant mean curvature.

Let $\mathcal{M}'_m$ and  $\mathcal{M}''_m$ be  meridian surfaces of elliptic and hyperbolic type, respectively.
Equality \eqref{E:Eq6} implies that the mean curvature of the meridian surface of elliptic type
$\mathcal{M}'_m$  is given by
\begin{equation} \label{E:Eq13}
|| H || = \sqrt{\frac{\varepsilon(\kappa^2 (\dot{f}^2 - 1) - (f \ddot{f} + \dot{f}^2 - 1)^2)}{4f^2 (\dot{f}^2 - 1)}}.
\end{equation}
Similarly, from \eqref{E:Eq10} it follows that the mean curvature of the meridian surface of hyperbolic type
$\mathcal{M}''_m$  is
\begin{equation} \label{E:Eq14}
|| H || = \sqrt{\frac{\varepsilon ((f \ddot{f} + \dot{f}^2 - 1)^2 - \kappa^2 (1 - \dot{f}^2))}{4f^2 (1 - \dot{f}^2)}}.
\end{equation}

\begin{thm}  \label{T:mean curvature}
Let $\mathcal{M}'_m$ (resp.  $\mathcal{M}''_m$) be a meridian surface of elliptic (resp. hyperbolic) type from the general class.

(i)  $\mathcal{M}'_m$  has constant mean curvature $|| H || = a = const$, $a \neq 0$ if and only if
the curve $c$  on $S^2(1)$ has
constant spherical curvature $\kappa = const = b, \; b \neq 0$,
and the meridian $m$ is determined by $\dot{f} = y(f)$
where
\begin{equation} \notag
y(t) =  \sqrt{1 + \frac{1}{t^2}\left(C \pm \frac{t}{2} \sqrt{b^2 - 4 a^2 t^2} \pm \frac{b^2}{4a}  \arcsin \frac{2at}{b}\right)^2}, \qquad C = const,
\end{equation}
$g(u)$ is defined by $\dot{g}(u) = \sqrt{\dot{f}^2(u) - 1}$.

(ii)  $\mathcal{M}''_m$  has constant mean curvature $|| H || = a = const$, $a \neq 0$ if and only if
 the curve $c$  on $S^2_1(1)$ has constant
spherical curvature $\kappa = const = b, \; b \neq 0$,
and the meridian $m$ is determined by $\dot{f} = y(f)$
where
\begin{equation} \notag
y(t) =  \sqrt{1 - \frac{1}{t^2}\left(C \pm \frac{t}{2} \sqrt{b^2 - 4 a^2 t^2} \pm \frac{b^2}{4a}  \arcsin \frac{2at}{b}\right)^2}, \qquad C = const,
\end{equation}
$g(u)$ is defined by $\dot{g}(u) = \sqrt{1 - \dot{f}^2(u)}$.
\end{thm}

\noindent {\it Proof:}
(i) It follows from \eqref{E:Eq13} that $||H|| = a$ if
and only if
$$\kappa^2(v) = \frac{4 a^2 f^2(u)(\dot{f}^2(u) - 1) + (f(u) \ddot{f}(u) + \dot{f}^2(u) - 1)^2}{\dot{f}^2(u) - 1},$$
which implies
\begin{equation} \label{E:Eq15}
\begin{array}{l}
\vspace{2mm}
\kappa = const = b, \; b \neq 0;\\
\vspace{2mm} 4 a^2 f^2(u)(\dot{f}^2(u) - 1) + (f(u) \ddot{f}(u) + \dot{f}^2(u) - 1)^2 = b^2(\dot{f}^2(u) - 1).
\end{array}
\end{equation}
The first equality of \eqref{E:Eq15} implies that the spherical curve $c$
has constant spherical curvature $\kappa = b$,
 i.e. $c$ is a circle on $S^2(1)$.
The second equality of \eqref{E:Eq15} gives the
following differential equation for the meridian $m$:
\begin{equation} \label{E:Eq16}
(f \ddot{f} + \dot{f}^2 - 1)^2 = (\dot{f}^2 - 1) (b^2 - 4 a^2 f^2).
\end{equation}

Further, if we set $\dot{f} = y(f)$ in equation \eqref{E:Eq16}, we obtain
that the function $y = y(t)$ is a solution of the following
differential equation
$$\frac{t}{2}(y^2)' + y^2 - 1 = \pm \sqrt{y^2 - 1} \sqrt{b^2 - 4 a^2 t^2}.$$
The general solution of the above equation is given by the formula
\begin{equation} \label{E:Eq17}
y(t) =  \sqrt{1 + \frac{1}{t^2}\left(C \pm \frac{t}{2} \sqrt{b^2 - 4 a^2 t^2} \pm \frac{b^2}{4a}  \arcsin \frac{2at}{b}\right)^2}, \qquad C = const.
\end{equation}
The function $f(u)$ is determined by $\dot{f} = y(f)$ and \eqref{E:Eq17}. The function $g(u)$ is
defined by $\dot{g}(u) = \sqrt{\dot{f}^2(u) - 1}$.

(ii) Similarly to the elliptic case, from \eqref{E:Eq14} it follows that  $||H|| = a$ if
and only if
the curve $c$ on $S^2_1(1)$ has constant curvature $\kappa = b$,
 and the meridian $m$ is determined by the following differential equation:
\begin{equation} \label{E:Eq18}
(f \ddot{f} + \dot{f}^2 - 1)^2  = (1 - \dot{f}^2) (b^2 + 4 a^2 f^2).
\end{equation}
Setting  $\dot{f} = y(f)$ in equation \eqref{E:Eq18}, we obtain
\begin{equation} \notag
y(t) =  \sqrt{1 - \frac{1}{t^2}\left(C \pm \frac{t}{2} \sqrt{b^2 - 4 a^2 t^2} \pm \frac{b^2}{4a}  \arcsin \frac{2at}{b}\right)^2}, \qquad C = const,
\end{equation}
In this case the function $g(u)$ is defined by $\dot{g}(u) = \sqrt{1 - \dot{f}^2(u)}$.

\qed

\section{Meridian surfaces with constant invariant $k$}

Let $\mathcal{M}'_m$ and  $\mathcal{M}''_m$ be  meridian surfaces of elliptic and hyperbolic type, respectively.
Then the invariant $k$ is given by the formula
\begin{equation} \label{E:Eq19}
k = - \frac{\kappa_m^2(u) \, \kappa^2(v)}{f^2(u)},
\end{equation}
where
$\kappa_m (u) = \ds{\frac{\ddot{f}(u)}{\sqrt{\dot{f}^2(u) - 1}}}$ in the elliptic case, and
$\kappa_m (u) = \ds{- \frac{\ddot{f}(u)}{\sqrt{1 - \dot{f}^2(u)}}}$ in the hyperbolic case.

In the following theorem we describe the meridian surfaces of elliptic or hyperbolic type with constant invariant $k$.

\begin{thm} \label{T:constant k}
Let $\mathcal{M}'_m$ (resp.  $\mathcal{M}''_m$) be a meridian surface of elliptic (resp. hyperbolic) type from the general class.

(i)  $\mathcal{M}'_m$  has constant invariant $k = const = - a^2, \; a\neq 0$ if and only if the curve $c$  on $S^2(1)$ has constant
spherical curvature $\kappa = const = b, \; b \neq 0$, and the meridian $m$ is determined by $\dot{f} = y(f)$
where
\begin{equation} \notag
y(t) =  \sqrt{1 + \left(C \pm \frac{at^2}{2b}\right)^2}, \qquad C = const,
\end{equation}
$g(u)$ is defined by $\dot{g}(u) = \sqrt{\dot{f}^2(u) - 1}$.

(ii)  $\mathcal{M}''_m$ has constant invariant $k = const = - a^2, \; a\neq 0$ if and only if the curve $c$  on $S^2_1(1)$ has constant
spherical curvature $\kappa = const = b, \; b \neq 0$,
and the meridian $m$ is determined by $\dot{f} = y(f)$
where
\begin{equation} \notag
y(t) =  \sqrt{1 - \left(C \mp \frac{at^2}{2b}\right)^2}, \qquad C = const,
\end{equation}
$g(u)$ is defined by $\dot{g}(u) = \sqrt{1 - \dot{f}^2(u)}$.
\end{thm}

\noindent {\it Proof:}
(i) It follows from \eqref{E:Eq19} that $k = const = -
a^2, \; a \neq 0$ if and only if
$$\kappa^2(v) = \frac{a^2 f^2(u) (\dot{f}^2(u) - 1)}{\ddot{f}\,^2(u)}.$$
The last equality implies
\begin{equation} \notag
\begin{array}{l}
\vspace{2mm}
\kappa = const = b, \; b \neq 0;\\
\vspace{2mm}
 a^2 f^2(u) (\dot{f}^2(u) - 1) = b^2 \ddot{f}\,^2(u).
\end{array}
\end{equation}
Hence, the curve $c$ has constant spherical curvature $\kappa = b$ and the function $f(u)$ is a solution of the following
differential equation:
\begin{equation} \label{E:Eq20}
b^2 \ddot{f}\,^2 - a^2 f^2 (\dot{f}^2 - 1) = 0
\end{equation}

Setting $\dot{f} = y(f)$ in equation \eqref{E:Eq20}, we obtain that
the function $y = y(t)$ is a solution of
\begin{equation} \notag
\frac{b}{2}(y^2)'  = \pm at \sqrt{y^2 - 1}.
\end{equation}
The general solution of the above equation is given by
\begin{equation}  \label{E:Eq21}
y(t) =  \sqrt{1 + \left(C \pm \frac{at^2}{2b}\right)^2}, \qquad C = const.
\end{equation}
The function $f(u)$ is determined by $\dot{f} = y(f)$ and \eqref{E:Eq21}.
The function $g(u)$ is defined by $\dot{g}(u) = \sqrt{\dot{f}^2(u) - 1}$.

(ii) Similarly to the elliptic case we obtain that  $\mathcal{M}''_m$ has constant invariant $k = const = - a^2, \; a\neq 0$ if and only if
 $c$  has constant curvature $\kappa = const = b, \; b \neq 0$,
and the meridian $m$ is determined by the following differential equation:
\begin{equation} \notag
b^2 \ddot{f}\,^2 - a^2 f^2 (1 - \dot{f}^2) = 0
\end{equation}

Again setting $\dot{f} = y(f)$  we obtain
\begin{equation} \notag
y(t) =  \sqrt{1 - \left(C \mp \frac{at^2}{2b}\right)^2}, \qquad C = const.
\end{equation}

\qed

\section{Chen meridian surfaces}

Let $\mathcal{M}'_m$ and  $\mathcal{M}''_m$ be  meridian surfaces of elliptic and hyperbolic type, respectively.
The invariants of $\mathcal{M}'_m$ and  $\mathcal{M}''_m$ are given by formulas \eqref{E:Eq7} and \eqref{E:Eq11}, respectively.
Recall that a spacelike surface in $\R^4_1$  is a non-trivial  Chen surface if and only if  $\lambda = 0$.
In the following theorem we classify all Chen meridian surfaces of elliptic or hyperbolic type.

\begin{thm} \label{T:Chen}
Let $\mathcal{M}'_m$ (resp.  $\mathcal{M}''_m$) be a meridian surface of elliptic (resp. hyperbolic) type from the general class.

(i)  $\mathcal{M}'_m$ is  a Chen surface  if and only if
the curve $c$  on $S^2(1)$ has
constant spherical curvature $\kappa = const = b, \; b \neq 0$,
and the meridian $m$ is determined by $\dot{f} = y(f)$
where
\begin{equation} \notag
y(t) =  \frac{\pm 1}{2\,t^{\pm1}} \sqrt{4\, t^{\pm2} - a \left(t^{\pm2} - \frac{b^2}{a}\right)^2}, \qquad a = const \neq 0,
\end{equation}
$g(u)$ is defined by $\dot{g}(u) = \sqrt{\dot{f}^2(u) - 1}$.

(ii)  $\mathcal{M}''_m$ is a Chen surface if and only if
 the curve $c$  on $S^2_1(1)$ has constant
spherical curvature $\kappa = const = b, \; b \neq 0$,
and the meridian $m$ is determined by $\dot{f} = y(f)$
where
\begin{equation} \notag
y(t) =  \frac{\pm 1}{2\,t^{\pm1}} \sqrt{4\, t^{\pm2} + a \left(t^{\pm2} - \frac{b^2}{a}\right)^2}, \qquad a = const \neq 0,
\end{equation}
$g(u)$ is defined by $\dot{g}(u) = \sqrt{1 - \dot{f}^2(u)}$.
\end{thm}

\noindent {\it Proof:}
(i) It follows from \eqref{E:Eq7} that $\lambda = 0$ if
and only if
\begin{equation} \notag
\kappa^2(v) =  \frac{ (\dot{f}^2(u) - 1)^2 - f^2(u) \ddot{f}\,^2(u)}{\dot{f}^2(u) - 1},
\end{equation}
which implies
\begin{equation} \notag
\begin{array}{l}
\vspace{2mm}
\kappa = const = b, \; b \neq 0;\\
\vspace{2mm}
(\dot{f}^2(u) - 1)^2 - f^2(u) \ddot{f}\,^2(u) = b^2 (\dot{f}^2(u) - 1).
\end{array}
\end{equation}

Hence, the curve $c$ has constant spherical curvature $\kappa = b$ and
the function $f(u)$ is a solution of the
following differential equation:
\begin{equation} \label{E:Eq22}
(\dot{f}^2 - 1)^2 - f^2 \ddot{f}\,^2 = b^2 (\dot{f}^2 - 1).
\end{equation}

The solutions of differential equation \eqref{E:Eq22} can be found as follows. Setting $\dot{f} = y(f)$ in equation \eqref{E:Eq22}, we obtain
that the function $y = y(t)$ is a solution of the equation:
\begin{equation} \label{E:Eq23}
\frac{t^2}{4} \left( (y^2)' \right)^2 = (y^2 - 1)^2 - b^2(y^2 - 1).
\end{equation}
We set $z(t) = y^2(t) - 1$ and obtain
\begin{equation} \notag
\frac{t}{2} \,z' = \pm \sqrt{z^2 - b^2 z}.
\end{equation}
The last equation is equivalent to
\begin{equation} \label{E:Eq24}
\frac{z'}{\sqrt{z^2 - b^2 z}} = \pm \frac{2}{t}.
\end{equation}
Integrating both sides of \eqref{E:Eq24}, we get
\begin{equation} \label{E:Eq25}
 \frac{b^2}{2} - z + \sqrt{z^2 - b^2 z} = c \, t^{\pm2}, \qquad c = const.
\end{equation}
It follows from \eqref{E:Eq25} that
\begin{equation} \notag
z(t) = - \frac{(a\, t^{\pm2} - b^2)^2}{4a \,t^{\pm2}}, \qquad a = 2c.
\end{equation}
 Hence, the general solution of differential equation \eqref{E:Eq23} is given by
\begin{equation} \notag
y(t) =  \frac{\pm 1}{2\,t^{\pm1}} \sqrt{4\, t^{\pm2} - a \left(t^{\pm2} - \frac{b^2}{a}\right)^2}, \qquad a = const \neq 0.
\end{equation}

(ii) In a similar way, in the hyperbolic case we obtain that  $\lambda = 0$ if and only if
the curve $c$ has constant curvature $\kappa = b$, $b \neq 0$ and
the function $f(u)$ is a solution of
\begin{equation} \notag
(1 - \dot{f}^2)^2 - f^2 \ddot{f}\,^2 = b^2 (1 - \dot{f}^2).
\end{equation}
Doing similar calculations as in the previous case, we obtain
\begin{equation} \notag
y(t) =  \frac{\pm 1}{2\,t^{\pm1}} \sqrt{4\, t^{\pm2} + a \left(t^{\pm2} - \frac{b^2}{a}\right)^2}, \qquad a = const \neq 0.
\end{equation}

\qed

\section{Meridian surfaces with parallel normal bundle}

In this section we shall describe the meridian surfaces of elliptic or hyperbolic type with parallel normal bundle.
Recall that a surface in $\R^4_1$ has parallel normal bundle if and only if $\beta_1 = \beta_2 =0$.

\begin{thm} \label{T:parallel}
Let $\mathcal{M}'_m$ (resp.  $\mathcal{M}''_m$) be a meridian
surface of elliptic (resp. hyperbolic) type from the general
class.

(i)  $\mathcal{M}'_m$ has parallel normal bundle if and only if
one of the following cases holds:

\hskip 10mm (a) the meridian $m$ is defined by
\begin{equation} \notag
\begin{array}{l}
\vspace{2mm}
f(u) = \pm \sqrt{u^2 + 2cu +d};\\
\vspace{2mm}
g(u) = \pm \sqrt{c^2 - d} \, \ln |u + c + \sqrt{u^2 + 2cu +d}| + a,
\end{array}
\end{equation}
where  $a$, $c$, and $d$ are constants, $c^2 > d$;

\hskip 10mm (b) the curve $c$   on $S^2(1)$ has constant spherical
curvature $\kappa = const = b, \; b \neq 0$, and the meridian $m$ is
determined by $\dot{f} = y(f)$ where
\begin{equation} \notag
y(t) = \pm \frac{\sqrt{(a^2+1)\,t^2 + 2ac\,t + c^2}}{t}, \quad a = const \neq 0, \quad c = const,
\end{equation}
$g(u)$ is
defined by $\dot{g}(u) = \sqrt{\dot{f}^2(u)-1}$.

(ii)  $\mathcal{M}''_m$ has parallel normal bundle if and only if
one of the following cases holds:

\hskip 10mm (a) the meridian $m$ is defined by
\begin{equation} \notag
\begin{array}{l}
\vspace{2mm}
f(u) = \pm \sqrt{u^2 + 2cu +d};\\
\vspace{2mm} g(u) = \pm \sqrt{d - c^2} \, \ln |u + c + \sqrt{u^2 +
2cu +d}| + a,
\end{array}
\end{equation}
where  $a$, $c$, and $d$ are constants, $d > c^2$;

\hskip 10mm (b) the curve $c$   on $S^2_1(1)$  has constant
spherical curvature $\kappa = const = b, \; b \neq 0$, and the
meridian $m$ is determined by $\dot{f} = y(f)$ where
\begin{equation} \notag
y(t) = \pm \frac{\sqrt{(1-a^2)\,t^2 + 2ac\,t - c^2}}{t},  \quad a = const \neq 0, \quad c = const,
\end{equation}
$g(u)$ is defined by $\dot{g}(u) = \sqrt{1 - \dot{f}^2(u)}$.
\end{thm}

\noindent {\it Proof:}
(i)  Using formulas \eqref{E:Eq7} we get
that $\beta_1 = \beta_2 =0$ if and only if
\begin{equation} \label{E:Eq26}
\begin{array}{l}
\vspace{2mm}
 \ds{\kappa\, \frac{d}{du}\left( \frac{f \ddot{f} + \dot{f}^2 - 1}{\sqrt{\dot{f}^2 - 1}} \right) -  \frac{d}{dv}(\kappa) \, \frac{f \ddot{f} + \dot{f}^2 - 1}{f \sqrt{\dot{f}^2 - 1}} = 0};\\
\vspace{2mm}
\ds{\kappa\, \frac{d}{du}\left( \frac{f \ddot{f} + \dot{f}^2 - 1}{\sqrt{\dot{f}^2 - 1}} \right) +  \frac{d}{dv}(\kappa) \, \frac{f \ddot{f} + \dot{f}^2 - 1}{f \sqrt{\dot{f}^2 - 1}} = 0}.
\end{array}
\end{equation}
It follows from \eqref{E:Eq26} that there are  two possible  cases:

\vskip 1mm
Case (a): $f \ddot{f} + \dot{f}^2 - 1 = 0$.
The general solution of this differential equation is  $f(u) = \pm \sqrt{u^2 + 2cu +d}$, $c = const$,  $d = const$.
Using that $\dot{g}^2 = \dot{f}^2-1$, we get  $\dot{g}^2 = \ds{ \frac{c^2 - d}{u^2 + 2cu +d}}$, and hence $c^2 - d >0$.
Integrating both sides of the equation
\begin{equation} \notag
\dot{g}(u) = \pm \ds{ \frac{\sqrt{c^2 - d}}{\sqrt{u^2 + 2cu +d}}},
\end{equation}
we obtain $g(u) = \pm \sqrt{c^2 - d}\, \ln |u + c + \sqrt{u^2 + 2cu +d}| + a$, $a = const$.
Consequently, the meridian $m$ is defined as described in \emph{(a)}.

\vskip 1mm
Case (b):
$\ds{ \frac{f \ddot{f} + \dot{f}^2 - 1}{\sqrt{\dot{f}^2 - 1}}} = a = const$, $a \neq 0$ and $\kappa = b = const$, $b \neq 0$.
Hence,  in this case the curve $c$ has constant spherical curvature $\kappa = b$ and
the meridian $m$ is determined by the following differential
equation:
\begin{equation} \label{E:Eq27}
f \ddot{f} + \dot{f}^2 - 1 = a \sqrt{\dot{f}^2 - 1}, \qquad a = const \neq 0.
\end{equation}

The solutions of differential equation \eqref{E:Eq27} can be found in the  following way. Setting $\dot{f} = y(f)$ in equation \eqref{E:Eq27}, we obtain
that the function $y = y(t)$ is a solution of the equation:
\begin{equation} \label{E:Eq28}
\frac{t}{2} \,(y^2)' + y^2 -1 = a \sqrt{y^2 - 1}.
\end{equation}
If we  set $z(t) = \sqrt{y^2(t) -1}$ we get
\begin{equation} \notag
z' + \frac{1}{t}\, z = \frac{a}{t}.
\end{equation}
The general solution of the above equation is given by the formula $z(t) = \ds{\frac{c + at}{t}}$, $ c = const$. Hence, the general solution of \eqref{E:Eq28} is
\begin{equation} \notag
y(t) = \pm \frac{\sqrt{(a^2+1)\,t^2 + 2ac\,t + c^2}}{t}, \quad c = const.
\end{equation}

(ii) In a similar way, considering meridian surfaces of hyperbolic type we obtain that  $\beta_1 = \beta_2 =0$  if and only if
one of the following cases holds.

\vskip 1mm
Case (a): $f \ddot{f} + \dot{f}^2 - 1 = 0$. In this case we get
\begin{equation} \notag
f(u) = \pm \sqrt{u^2 + 2cu +d}; \qquad
g(u) = \pm \sqrt{d - c^2} \, \ln |u + c + \sqrt{u^2 + 2cu +d}| + a,
\end{equation}
where  $a$, $c$, and $d$ are constants, $d > c^2$.

\vskip 1mm
Case (b):
$\ds{ \frac{f \ddot{f} + \dot{f}^2 - 1}{\sqrt{1 - \dot{f}^2}}} = a = const$, $a \neq 0$ and $\kappa = b = const$, $b \neq 0$.
Doing similar calculations as the calculations for solving \eqref{E:Eq27}, we obtain
\begin{equation} \notag
y(t) = \pm \frac{\sqrt{(1-a^2)\,t^2 + 2ac\,t - c^2}}{t}, \quad c =
const.
\end{equation}

\qed

\vskip 3mm
Similarly to the elliptic or hyperbolic type one can study the invariants of the meridian surfaces of parabolic type.
The classes of meridian surfaces of parabolic type  with constant Gauss curvature, constant mean curvature,
constant invariant $k$, the Chen meridian surfaces of parabolic type,
and the meridian surfaces of parabolic type with parallel normal bundle can be described in an analogous way.

\end{document}